# Theorems on Separated Transformations of Basis Vectors of Polynomial Space and its Applications in Classical Orthogonal Polynomials

Manouchehr Amiri[1]
Independent Researcher

The present paper introduces a method of basis transformation of a vector space that is specifically applicable to polynomials space and differential equations with certain polynomials solutions such as Hermite, Laguerre and Legendre polynomials. This method is based on separated transformations of vector space basis by a set of operators that their overall action is equivalent to the formal basis transformation and connected to it by linear combination with projection operators. By introducing the new relations for projection operators under separated transformation and applying the Forbenius covariants as projection operators, we derive a general method that incorporates the Rodrigues formula as a special case in polynomial space. A solution method to some classical differential equation such as Hermite and Laguerre differential equations is presented.

## 1. Introduction

In general mathematics, the solution of many problems requires solving the differential equations and their eigenvalues problem . Hermite, Laguerre and Legendre polynomials and related differential equations are among the most applicable eigenvalues problem in physics and mathematics [1, 2, 6]. Schrodinger equation for hydrogen atom reduces to Legendre differential equation and quantum harmonic oscillator requires Hermite polynomials and related differential equation. The well-known Rodrigues formula yields the solutions (eigenfunctions) of many of these differential equations [3]. In this paper we interpret the Rodrigues formula as the transformation of some specific basis in polynomial space to another polynomials (eigenfunctions) regarding the associated differential equations. This approach is feasible, provided that the transformation operator to be considered as a set of operators acting on each basis separately. It is shown that the overall action of these operators is equivalent to a single linear operator. As an example, the change of basis vectors in two dimension can be made by a matrix of rank 2. The action of this matrix could be equivalent with the actions of two different matrices that acts on each basis separately. through a theorem in section 2 the general form of projection operators under the separated transformation operators is proved. The relation between these operators could be achieved by applying projection operators as is also proved in the section 2. The connection between separated basis transformation and umbral composition has been revealed by a theorem in section 2. In section 3, it is proved that by knowing the first two polynomials of Hermite , Laguerre and Legendre polynomials, the related differential equations could be retrieved  using the method based on the separated transformation of original basis and Forbenius covariants [4, 5] as projection operators. The examples in section 3 clarifies the details of this method. By applying the Rodrigues formula as the separated operators acting on the

manoamiri@gmail.com [1]

original basis, we acquire the form of related differential equations. In section 3.3, we propose a technique for solutions of differential equations based on raising operators acquired from associated Lie algebras.

## 2. Transformation of Basis Vectors by Separated Operators

Let $\mathbb{V}$ be a n-dimensional vector space with basis vectors $e_1, e_2, \ldots e_n$. The linear operator that transforms these basis to another basis $e'_1, e'_2, \ldots e'_n$, normally is defined as a unique linear operator $O$ in the matrix form. In present theory we define a set of linear operators $O_1, O_2, \ldots O_n$; each one acts **separately** on the corresponding basis as follows:
$$O_i e_i = e'_i \tag{1}$$
The result of the action of operator $O$ and the set of $O_i$ on the initial basis are the same, but $O_i$ as separated basis transformations allow to choose a wide range of $O_i$ operators whose overall transformations are equivalent to the operation of $O$. On the other hand, in many problems such as Rodrigues type formulas the separated basis transformation $O_i$ are more accessible than overall operator $O$. In the context of differential operators, $O_i$ could be regarded as the operators that transform the initial basis (monomials) in polynomial space to another basis as for example we observe in Rodrigues formula for Laguerre polynomials as the solutions (eigenfunctions) to Laguerre differential equation [6]:
$$\mathbb{L}_n = \frac{1}{n!}(D-1)^n x^n$$
This equation can be interpreted as the transformation of initial basis $(1, x, x^2, x^3, \ldots)$ to new basis $\mathbb{L}_n$ (Laguerre polynomials) by the action of the operator,
$$O_n = \frac{1}{n!}(D-1)^n$$
Respect to Rodrigues formula, for all related differential equations such as Legendre, Chebyshev and Bessel Equations, there are separate and independent operators for each basis. Therefor we can apply a set of operators $O_i$ instead of a unique operator $O$ to transform the initial basis in polynomials space. This method obviates the need to find the unique linear operator $O$ with the same action on all initial bases. We will show the relation between $O_i$ and $O$ in proposition 2.1 after defining the projection operators as follows.
Let introduce the projection operators $P_1, P_2, \ldots, P_n$ by the definition:

$$P_i V = V_i e_i \; ; \; P_i e_j = \delta_{ij} e_i \tag{2}$$

Where $V \in \mathbb{V}$ is a vector expanded as:
$$V = \sum_i V_i e_i \tag{3}$$
From (2) we have: $\qquad \sum_i P_i = I \; ; \; P_i P_j = P_i \delta_{ij}$
as the main condition for projection operators ($I$ is identity operator).

**Remark 1**
projection operators defined in (2), are linear operators.
We show that by basis transformations according to (1), the projection operators $P_i$ are transformed as:

$$P'_i = O_i P_i (\sum_j O_j P_j)^{-1} \qquad (4)$$

**Theorem 2.1**: the generalized form of projection operator under **separated** basis transformations $O_i e_i = e'_i$ is:

$$P'_i = O_i P_i (\sum_j O_j P_j)^{-1}$$

**Proof**: respect to basis transformation $O_i e_i = e'_i$ and (2) we have:

$$e'_i = P'_i e'_i = P'_i O_i e_i = P'_i (\sum_j O_j P_j) e_i \qquad (5)$$

Again, with substitution $e'_i = O_i e_i$ (5) reads as:

$$O_i e_i = P'_i (\sum_j O_j P_j) e_i \qquad (6)$$

It is valid for all $e_i$. with identity $O_i e_i = O_i P_i e_i$, equation (6) becomes:

$$O_i P_i = P'_i (\sum_j O_j P_j) \qquad (7)$$

Or:
$$P'_i = O_i P_i (\sum_j O_j P_j)^{-1} \qquad \square$$

As we expected.

**Proposition 2.1** The transformation of projection operator defined in (4) is equivalent to a similarity transformation:

$$P'_i = (\sum_k O_k P_k) P_i (\sum_j O_j P_j)^{-1} = O_i P_i (\sum_j O_j P_j)^{-1} \qquad (8)$$

**Proof**: expansion of the first two terms on left side considering $P_i P_j = P_i \delta_{ij}$ yields:

$$(\sum_k O_k P_k) P_i = O_i P_i \qquad (8a)$$

Therefor we have the right side of (8). $\square$

Equation (8) implies the similarity transformation of $P_i$ under the basis transformation made by operator i.e.

$$O = \sum_i O_i P_i . \qquad (8b)$$

Therefor the operator $O$ is the linear operator for transforming all basis $e_i$. Operator $O$ also yields transformations of all linear operators $\mathcal{O}$ in vector space $\mathbb{V}$ under basis transformation $e'_i = O e_i$ by the similarity transformation

$$\mathcal{O}' = O \mathcal{O} O^{-1} = (\sum_k O_k P_k) \mathcal{O} (\sum_j O_j P_j)^{-1}$$

**Remark 2**

Equation (4) meets the projection operator conditions:

a) $\sum_i P'_i = \sum_i O_i P_i (\sum_j O_j P_j)^{-1} = (\sum_i O_i P_i)(\sum_j O_j P_j)^{-1} = I \qquad (9)$

$$P'_i = O P_i O^{-1} = \sum_k O_k P_k) P_i (\sum_j O_j P_j)^{-1} = O_i P_i (\sum_j O_j P_j)^{-1}$$

The action of $O$ and $O_j$ on a basis $e_j$ is the same.

$$O e_j = (\sum_i O_i P_i) e_j = O_j P_j e_j = O_j e_j = e'_i$$

This implies that the action of $O$ and $O_j$ on $e_j$ is equivalent.

b) Respect to (4) and (9) we conclude:
$$P'_i P'_j = O_i P_i (\sum_k O_k P_k)^{-1} O_j P_j (\sum_l O_l P_l)^{-1} =$$
$$[(\sum_k O_k P_k) P_i (\sum_j O_j P_j)^{-1}][(\sum_k O_k P_k) P_j (\sum_j O_j P_j)^{-1}] \qquad (10)$$

Middle terms in right side of (7) reduce to identity operator and thus:
$$P'_i P'_j = (\sum_k O_k P_k) P_i P_j (\sum_j O_j P_j)^{-1} \qquad (11)$$

Recalling $P_i P_j = P_i \delta_{ij}$ and (4) we obtain:
$$P'_i P'_j = O_i P_i (\sum_j O_j P_j)^{-1} \delta_{ij} = P'_i \delta_{ij}$$

This proves the *idempotency* of $P'_i$ i.e., $\quad P'_i P'_i = P'_i$

Where $P'_i$ denoted as posterior probability analogy.

Equation (4) is the unique formula for $P'_i$ and other forms in spite of their validity for satisfaction of projection operator conditions i.e., equation (3), are not the right candidates. As an example, we may propose this formula for $P'_i$:

$$P'_i = (\sum_j P_j O_j)^{-1} P_i O_i \tag{12}$$

It is straightforward to investigate that this definition is compatible with conditions (3) but if we multiply both sides by $e'_i$ we obtain:

$$(\sum_j P_j O_j) P'_i e'_i = P_i O_i e'_i$$

Respect to (1) and (2) we get:

$$(\sum_j P_j O_j) e'_i = P_i O_i e'_i$$

One of the solutions results in a false outcome:

$$(\sum_j P_j O_j) = P_i O_i$$

**Proposition 2.2** The product of projection operators is associative.

**Proof** : Let projection operators $P'_i$ and $P''_i$ correspond to $O_i$ and $O'_i$ as transformation groups of coordinates i.e.

$$P'_i = O_i P_i (\sum_j O_j P_j)^{-1}$$

And
$$P''_i = O'_i P'_i (\sum_j O'_j P'_i)^{-1}$$

Substitution of first relation into above equation results in:

$$P''_i = O'_i O_i P_i (\sum_j O_j P_j)^{-1} \; [(\sum_j O'_j \, O_j P_j)(\sum_j O_j P_j)^{-1})]^{-1}$$

$$P''_i = O'_i O_i P_i (\sum_j O_j P_j)^{-1} [(\sum_j O_j P_j)(\sum_j O'_j \, O_j P_j)^{-1}]$$

After vanishing of two central terms:

$$P''_i = O'_i O_i P_i (\sum_j O'_j \, O_j P_j)^{-1}$$

This is compatible with (4) by replacing $O_i$ with $O''_i = O'_i O_i$. Therefore, the corresponding projection operator for two consecutive transformation $O_i$ and $O'_i$ is equivalent the projection operator of $O''_i = O'_i O_i$ transformation.

$$P''_i = O''_i P_i (\sum_j O''_i P_j)^{-1} \qquad \square$$

**Remark 3**

As is proved, the operators $\sum_j O_j P_j$ and $O$ are equivalent operators. Respect to equation (9), the projection operator $P_i$ transforms as a similarity transformation under the action of operator $\sum_j O_j P_j$, therefor the initial basis should be transformed by this operator and consequently $\sum_j O_j P_j$ and $O$ are equivalent.

We show the identity (4) is also valid in function space, where the linear projection operators are defined.

**Proposition 2.3** Let $\mathbb{V}$ be a n-dimensional function space over the real field $F$ with a set of orthogonal basis functions $\varphi_i$ and inner product defined on a closed interval $[a, b]$. The same definition in (2) can be applied on these basis

$$P_i F = c_i \varphi_i = \varphi_i \int_a^b \varphi_i F \, dx \tag{13}$$

Where
$$c_i = \langle \varphi_i, F \rangle = \int_a^b \varphi_i F \, dx \tag{14}$$

Are the coefficients in expansion of square integrable function $F$ in the basis $\varphi_i$ calculated by inner product of $\varphi_i$ and $F$ over the interval $[a, b]$.

**Proof**: with the identity (4) we conclude:

Then we have:
$$P'_i(\sum_j O_j P_j) = O_i P_i$$
$$P'_i(\sum_j O_j P_j)F = O_i P_i F$$
Respect to (13) and (14) we have:
$$P'_i(\sum_j O_j \varphi_j \int \varphi_j F \, dx) = O_i \varphi_i \int_a^b \varphi_i F \, dx \tag{15}$$
Regarding equation (1) we can choose $\varphi'_i = O_i \varphi_i$ as transformed basis that results in:
$$P'_i(\sum_j \varphi'_j \int_a^b \varphi_j F \, dx) = \varphi'_i \int_a^b \varphi_i F \, dx \tag{16}$$
With the definition in equation (2) for projection operator we have:
$$P'_i \varphi'_j = \delta_{ij} \varphi'_j$$
Therefor the equation (16) respect to (14) reads as:
$$\varphi'_i \int_a^b \varphi_i F \, dx = c_i \varphi'_i$$
$$c_i \varphi'_i = c_i \varphi'_i$$
So, the identity (9) is valid for operators in function spaces. □

**Proposition 2.4** Differential operators with certain eigenvalues and eigenfunctions can be linearly expanded by their projection operators.

**Proof:** Let the differential operator $\mathfrak{D}$ is characterized by eigenfunctions relation:
$$\mathfrak{D}\varphi_i = \lambda_i \varphi_i \tag{17}$$
The eigenfunctions $\varphi_i$ are linearly independent and are the basis vectors, i.e.: $P_i \varphi_j = \delta_{ij} \varphi_j$.
Where $P_i$ is the projection on $i$-th subspace, then by the identity:
$$\mathfrak{D}\varphi_i = \lambda_i \varphi_i = (\sum_j \lambda_j P_j) \varphi_i \tag{18}$$
The validity of this equation for all $\varphi_i$ yields:
$$\mathfrak{D} = \sum_j \lambda_j P_j \tag{19}$$
That proves the proposition. □

**Theorem 2.5** Let the initial basis $e_i$ correspond to some set of linearly independent non-homogenous polynomials such as the regular bases $(1, x, x^2, x^3, \dots)$. After transforming the bases by equation $O_i e_i = e'_i$ to new bases $e'_i$ which correspond the new linearly independent polynomials $P_n(x)$, if $\mathfrak{D}$ denoted as the differential operator with $e_i$ or equivalently $x^n$ (n-th exponent of $x$) as its eigenfunctions (or eigenvector), then the corresponding differential operator $\mathfrak{D}'$ with eigenfunctions $P_n(x)$ can be obtained by the relation:
$$\mathfrak{D}' = (\sum_k \lambda'_k O_k P_k)(\sum_j O_j P_j)^{-1} \tag{20}$$
Where $\lambda'_k$ are eigenvalues of $\mathfrak{D}'$.

**Proof** Respect to equation (19) the expansion of $\mathfrak{D}'$ in terms of $P_i'$ reads as:
$$\mathfrak{D}' = \sum_i \lambda'_i P_i' \tag{21}$$
Where $P_j'$ are projection operators onto the i-th subspace (i.e., $e'_i$). Substitution of $P_i'$ in equation (21) by equation (4) results in:
$$\mathfrak{D}' = \sum_i \lambda'_i O_i P_i (\sum_j O_j P_j)^{-1} = (\sum_i \lambda'_i O_i P_i)(\sum_j O_j P_j)^{-1} \tag{22}$$
This proves the theorem. □

## 2.1 Projection operators in terms of resolvents:

Associated to any differential operator in Hilbert space there are projection operators in terms of their resolvent i.e.
$$P_i = \int_{C_{v_i}} \frac{d\lambda}{(\lambda I - \mathcal{D})^{-1}} \quad \text{and} \quad P'_i = \int_{C_{v'_i}} \frac{d\lambda}{(\lambda I - \mathcal{D}')^{-1}} \tag{23}$$

Using (21), (22) we obtain:
$$\int_{C_{v'_i}} \frac{d\mu}{\mu I - \mathcal{D}'} = (O_i \int_{C_{v_i}} \frac{d\lambda}{\lambda I - \mathcal{D}})(\sum_j O_j \int_{C_{v_j}} \frac{d\lambda}{\lambda I - \mathcal{D}})^{-1} \tag{24}$$

For a unique transformation $O = O_i$ for all $\varphi_i$ we get:
$$\int_{C_{v'_i}} \frac{d\mu}{\mu I - \mathcal{D}'} = O(\int_{C_{v_i}} \frac{d\lambda}{\lambda I - \mathcal{D}}) O^{-1}$$

with expansion of resolvents $\frac{d\lambda}{\lambda I - \mathcal{D}}$ as a Neumann infinite series (polynomial), it is proved that the corresponding differential operator after action of operator $O$ on the base functions $\varphi_i$ as defined in proposition 2.1, can be presented by a similarity transformation:
$$\mathcal{D}' = O \mathcal{D} O^{-1} \tag{25}$$

**Example 2.1**: Eigenfunctions of the differential operator $\mathcal{D} = \frac{d}{dx}$ could be found as $\varphi_n = e^{nx}$. Transforming by $\varphi'_n = O \varphi_n = x \varphi_n = x e^{nx}$, The resulting corresponding differential operator respect to proposition 2.1 after substituting $O = x$ reads as:
$$\mathcal{D}' = x \mathcal{D} x^{-1} = x \left( \frac{-1}{x^2} + x^{-1} \mathcal{D} \right) = \frac{-1}{x} + \mathcal{D} \tag{26}$$

Action of this operator on $xe^{nx}$ gives:
$$(\frac{-1}{x} + \mathcal{D}) x e^{nx} = -e^{nx} + e^{nx} + nx e^{nx} = nx e^{nx} \tag{27}$$

Thus, the eigenfunctions of this operator are $xe^{nx}$ as expected. Because of the similarity relations of operators $\mathcal{D}$ and $\mathcal{D}'$, their eigenvalues are identical.

**Theorem 2.6** Let the linearly independent monomials $p_m(x)$ and $q_m(x)$ of polynomials $P_n(x)$ and $Q_n(x)$ of degree $n$ are connected by the operators $O_i$ as defined in equation (1) i.e.,
$$q_m(x) = O_m p_m(x) \tag{28}$$
Denote $P_m$ as projection operators that project functions of variable $x$ on basis $p_m(x)$ with the definition of equation (2) $\quad P_m p_n(x) = \delta_{mn} p_n(x)$
Then the operator $O = \sum_m O_m P_m$ acts as *umbral* composition on polynomial $P_n(x)$.

**Proof:** Let expand $P_n(x)$ in terms of monomial basis $p_m(x)$
$$P_n(x) = \sum_n a_{nm} p_m(x) \tag{29}$$
Then action of $O$ on $P_n(x)$ gives
$$O P_n(x) = \sum_i O_i P_i (\sum_m a_{nm} p_m(x)) = \sum_i O_i a_{ni} p_i(x) = \sum_i a_{ni} q_i(x) \tag{30}$$

This implies that the action of $O$ on $P_n(x)$ replaces the monomials $p_m(x)$ with $q_i(x)$ while the coefficients $a_{nm}$ in the expansion remains unchanged. This means that the operation of $O$ is equivalent with umbral composition by the definition:
$$P_n \circ Q = \sum_m a_{nm} q_m(x) \tag{31}$$
This definition coincides the action of $O$ on $P_n(x)$. $\square$

### 2.2 Forbenius covariant of operators

For other representation of projection operator in terms of differential operator we apply the Forbenius covariants [5] as projection operators (matrices) which are the coefficient of

*Sylvester's* formula. For a differential operator $\mathfrak{D}$ in polynomial space, the projection operator on the one-dimensional eigenfunction subspaces is given by

$$P_l = \prod_{k=1}^{n} \frac{\mathfrak{D}-\lambda_k}{\lambda_l-\lambda_k} \qquad k \neq l \tag{32}$$

These operators act on the functions in function space and yields their projections on basis $\varphi_i$ which are the eigenfunctions of $\mathfrak{D}$ with corresponding eigenvalues $\lambda_i$.

## 2.3 Similarity transformation:

Respect to equation (9) and related proposition, if we substitute $O_i P_i$ with $(\sum_k O_k P_k)$ respect to the identity, $P_i P_j = \delta_{ij} P_i$ we have:

$$P'_i = (\sum_k O_k P_k) P_i (\sum_j O_j P_j)^{-1} = O_i P_i (\sum_j O_j P_j)^{-1} \tag{33}$$

This equation is a similarity transformation of $P_i$ under the operator $\sum_k O_k P_k$. This similarity transformation corresponds to the basis transformation $O_i e_i = e'_i$. Actually, $\sum_k O_k P_k$ as an operator $\hat{O}$ transforms all basis $e_i$ to $e'_i$ and corresponds the coordinate transformation matrix. From this equation we can deduce similarity transformation for other operators provided that the operators in similarity transformation have common eigenvalues. Therefor the differential operators with *identical eigenvalues* could be related by similarity transformation. As an example, differential operator $\mathfrak{D} = x\frac{d}{dx} = xD$ with basis (eigenfunction) $\varphi_n = x^n$ transforms to another differential operator $\mathfrak{D}'$ with eigenfunction $\varphi'_i$ after the basis transformation $\varphi'_i = O_i \varphi_i$. Therefor we have the similarity transformation:

$$\mathfrak{D}' = (\sum_k O_k P_k) \mathfrak{D} (\sum_j O_j P_j)^{-1} \tag{34}$$

If all $O_j$ are the same namely $O$, (34) will be reduced to:

$$\mathfrak{D}' = (O \sum_k P_k) \mathfrak{D} (O \sum_j P_j)^{-1}$$
$$\mathfrak{D}' = O \mathfrak{D} O^{-1}$$

In these cases that the single operator transforms all bases, the exact closed form of related differential operator could be derived by this method. However, for cases with separate $O_i$, the validity of the retrieved differential operator relies on the action on the first two polynomials as we show in next sections. The following example clarifies the method.

**Example 2.2**:

Let the vector space $\mathbb{V}$ spanned by the linearly independent basis $(1, e^x, e^{2x}, \ldots)$ which are the eigenfunctions of operator $\mathfrak{D} = \frac{d}{dx}$. If these basis transforms to the new set of basis by multiplying with $e^{\frac{x^2}{2}}$ i.e., $(e^{\frac{x^2}{2}}, e^{\frac{x^2}{2}+x}, e^{\frac{x^2}{2}+2x}, \ldots)$ then the corresponding operator with these new basis as its eigenfunctions could be obtained by (34). In this case $O_k = O = e^{\frac{x^2}{2}}$. Thus, the equation (34) reduces to:

$$\mathfrak{D}' = O \mathfrak{D} O^{-1}$$
$$\mathfrak{D}' = e^{\frac{x^2}{2}} \mathfrak{D} e^{\frac{-x^2}{2}}$$

The term $\mathfrak{D} e^{\frac{-x^2}{2}}$ is not just the derivative of $e^{\frac{-x^2}{2}}$, but an operator that is equal to:

$$\mathfrak{D} e^{\frac{-x^2}{2}} = \frac{d}{dx}(e^{\frac{-x^2}{2}}) + e^{\frac{-x^2}{2}} \mathfrak{D}$$

Then we have:
$$\mathcal{D}' = e^{\frac{x^2}{2}} \mathcal{D} e^{\frac{-x^2}{2}} = e^{\frac{x^2}{2}}[\frac{d}{dx}(e^{\frac{-x^2}{2}}) + e^{\frac{-x^2}{2}}\mathcal{D}]$$
$$\mathcal{D}' = e^{\frac{x^2}{2}}(-xe^{\frac{-x^2}{2}} + e^{\frac{-x^2}{2}}\mathcal{D})$$
$$\mathcal{D}' = (-x + \mathcal{D})$$

The eigenfunctions of this operator are $e^{\frac{x^2}{2}+nx}$ with eigenvalues $n$ as expected. It is noteworthy to note that the expression for probabilist's Hermite polynomial $H_{e1}$ with the definition:
$$H_{e1} = e^{\frac{x^2}{2}} \frac{d}{dx} e^{\frac{-x^2}{2}}$$

Differs from $\mathcal{D}'$, because in this definition the term $\frac{d}{dx} e^{\frac{-x^2}{2}}$ is not an operator but merely the derivative of $e^{\frac{-x^2}{2}}$.

It is easy to prove that any function of $\mathcal{D}'$ can be expanded in terms of $\mathcal{D}$ as follows:

$$f(\mathcal{D}') = O f(\mathcal{D}) O^{-1}$$

## 3  Separated Basis Transformation Method based on Forbenius Covariants and its applications

Another approach to find $\mathcal{D}'$ in terms of $\mathcal{D}$ and $O_j$ is to apply the Forbenius covariant operators as projection operators mentioned in (32).
$$P_k = \prod_{l=1}^{N} \frac{\mathcal{D}-\lambda_l}{\lambda_k-\lambda_l} \qquad l \neq k \tag{35}$$

These operators are projection operators onto the k-th one-dimensional sub-space (basis) [8]. substituting these projectors in equation (34) results in

$$\mathcal{D}' = \left(\sum_k O_k \prod_{l=1}^{N} \frac{\mathcal{D}-\lambda_l}{\lambda_k-\lambda_l}\right) \mathcal{D} (\sum_j O_j \prod_{l=1}^{N} \frac{\mathcal{D}-\lambda_l}{\lambda_j-\lambda_l})^{-1} \qquad l \neq j \tag{36}$$

N denoted as the dimension of function or polynomial space.
This method in comparison with previous methods are more applicable because the calculation of inverse of a product of differential operators is easier than other methods.

In the following sections we introduce an applicable method to find $\mathcal{D}'$ in terms of $\mathcal{D}$. Taking into account the equation (22) we have:
$$\mathcal{D}' = (\sum_i \lambda'_i O_i P_i)(\sum_j O_j P_j)^{-1} \tag{37}$$

Replacing $P_i$ in this equation with Forbenius covariant operators in equation (35) as projection operators, yields a useful formula for $\mathcal{D}'$ as follows:

$$\boxed{\mathcal{D}' = (\sum_i \lambda'_i O_i \prod_{l=1}^{N} \frac{\mathcal{D}-\lambda_l}{\lambda_i-\lambda_l})(\sum_j O_j \prod_{l=1}^{N} \frac{\mathcal{D}-\lambda_l}{\lambda_j-\lambda_l})^{-1}}$$
(38)

This equation is applicable in next sections.
If all $O_i$ are the same i.e., $O_j = O$, we have:
$$\mathcal{D}' = O(\sum_i \lambda'_i P_i) O^{-1} \tag{39}$$

The condition of identical eigenvalues for $\mathcal{D}'$ and $\mathcal{D}$ is not necessary in equation (37) and the case of identical eigenvalues are special case of this equation. we apply this equation restricted to the first two polynomials i.e., two-dimensional polynomial space. Substitution of $P_i$ in (37) by

Forbenius covariants (35) yields an applicable method as we will show in examples. It is noteworthy to recall that the term $\sum_j O_j P_j$ stands for a linear operator (equivalent to a matrix) that transforms the basis $(1, x, x^2, \ldots)$ of polynomials space to another basis. For example, it transforms basis $(1, x, x^2, \ldots)$ to Hermite polynomial $H_{en}$ as new linearly independent basis by the techniques that is presented in next section.

## 3.1 Applications of Separated Basis Transformation Method in Polynomial Space

In the Sturm Liouville problem and related differential equations and their specific solutions such as Hermite, Laguerre, Legendre and Jacobi polynomials, the transformation of basis in function space seems to be an interesting subject. For example, transformation of basis $(1, x, x^2, x^3, \ldots)$ under the multiple derivatives which is compatible with Rodrigues' formula to derive Hermite polynomials, presented as follows [6]:

$$H_{en}(x) = e^{\frac{-D^2}{2}} x^n \tag{40}$$

Where $D = \frac{d}{dx}$ and $H_{en}$ are probabilistic Hermite polynomials.

In the case of Laguerre polynomials, we have the transformation:

$$L_n(x) = \frac{1}{n!}(D-1)^n x^n \tag{41}$$

Respect to our theory, these transformations are compatible with the operator action of $O_n$ separately on basis $(1, x, x^2, x^3, \ldots)$, for Hermite polynomial we have:

$$O_n = O = e^{\frac{-D^2}{2}} \tag{42}$$

And for Laguerre polynomials:

$$O_n = \frac{1}{n!}(D-1)^n \tag{43}$$

We introduce the operator $xD$ as the unique operator with basis $(1, x, x^2, x^3, \ldots)$ as its associated eigenfunctions with eigenvalues $(0,1,2,\ldots)$:

$$xD(x^n) = nx^n \tag{44}$$

Therefor we can use the equation (38) to find the differential operator that its polynomials are determined by applying related $O_n$ on basis $(1, x, x^2, x^3, \ldots)$ as in (40) and (41). By substitution of $P_k$ and $O_n$ in equations (35), (37) and (38) and using $\mathfrak{D} = xD$ in equation (38) we recover the corresponding differential equations with eigenfunctions such as $H_{en}$ and $L_n$ as presented in next examples. The presented technique uses the first two polynomials i.e., the two-dimensional space of polynomials with monomials of order 1 and 0. This facilitates the calculation of desired differential equations and shows that if an infinite set of polynomials present the eigenfunctions of a unique differential operator, then applying this technique for the first two eigenfunctions gives the exact form of related differential equation that holds true for all other eigenfunctions. We clarify this method by the following proposition.

**Proposition 3.1** Let the function $B(x)$ be an arbitrary differentiable function and $B^n(x) = e_n$ be the $n$-th integer power of $B(x)$. Assume the set of original basis $(1, B(x), B^2(x), \ldots B^n(x))$ for polynomial space, are the eigenfunctions of differential operator $\mathfrak{D}$ and these basis change to the set of linearly independent polynomials basis $\mathbb{P}_n = e'_n$ after separated transformations under the action of operators $O_n$ as defined in equation (1):

$$\mathbb{P}_n = O_n B^n(x) \qquad n \in \mathbb{N} \cup \{0\} \tag{45}$$

If the $\mathbb{P}_n$ are the eigenfunctions of the *n*-th order differential operator $\mathfrak{D}'$, of the form:
$$\mathfrak{D}'y = (C_2 D^2 + C_1 D)y = \beta y \tag{46}$$

Where $C_m$ are the polynomials of order $m \leq 2$ and $D = \frac{d}{dx}$, then applying the Forbenius covariant operator defined in (35) and equation (38) for the first 2 polynomials (eigenfunctions) $\mathbb{P}_0$ and $\mathbb{P}_1$ is sufficient to yield the corresponding differential operator $\mathfrak{D}'$ derived from $O_n$ and $P_n$ as was derived in equation (38).

**Proof:** Acting the differential operator $\mathfrak{D}'$ in equation (38) on the first 2 polynomials; $\mathbb{P}_0$ and $\mathbb{P}_1$ is sufficient for determining the 2 coefficients $C_m$ because the last term in operator $\mathfrak{D}'$ in equation (38) is equivalent to the operator $O^{-1}$ as the inverse of $O$, defined in equation (8b):
$$(\textstyle\sum_j O_j \prod_{l=1}^N \frac{\mathfrak{D}-\lambda_l}{\lambda_j - \lambda_l})^{-1} = O^{-1} = (\textstyle\sum_i O_i P_i)^{-1}$$
and its action on $\mathbb{P}_0$ and $\mathbb{P}_1$ results in $1, B(x)$. Thus the 2 resulting equations are sufficient to determine the from of $\mathfrak{D}'$.

As a special case, we prove the Rodrigues formula as the action of operators $O_n$ on the initial basis $B^n(x)$ in polynomial space to transform them to new basis $\mathbb{P}_n$ that correspond to the desired differential operator $\mathfrak{D}'$(i.e., differential equation) as its eigenfunctions.

**Remark 4**

Applying this method for any subspace of the polynomial space holds true. As was shown in proposition 3.1, the differential operator $\mathfrak{D}'$ that is valid for 2-dimensional space spanned by $\mathbb{P}_0$ and $\mathbb{P}_1$, also valid for any basis $\mathbb{P}_n$ and therefore this method holds true for any subspace spanned by $(\mathbb{P}_0, \mathbb{P}_1, \ldots, \mathbb{P}_m)$; $m \leq n$. The example 3.1 clarifies this claim for 2- and 3-dimensional subspace.

### 3.2 Rodrigues' formula as a special case of separated basis transformation

In this section we prove the compatibility of Rodrigues' formula with our presented techniques and show that substitution of $O_n$ in equation (45) by Rodrigues' formula transformation, yields the corresponding differential operators and equations.

**Proposition 3.2** The corresponding differential operator $\mathfrak{D}'$ for transformation of basis vectors $B^n(x)$ defined by Rodrigues' formula as separated basis transformation:
$$\mathbb{P}_n(x) = \frac{1}{\omega} D^n[\omega B^n(x)]$$
is presented in the form:
$$\mathfrak{D}' = BD^2 + AD$$
Where $A$ is a polynomials of degree 1.

**Proof:** Due to the presented theory, we showed that if the bases $e_n$ of a vector space $\mathbb{V}$ which are the eigenfunctions of differential operator $\mathfrak{D}$, are transformed separately by operators $O_n$, the transformed differential operator obeys the equation (37) i.e.,

$$\mathfrak{D}' = (\textstyle\sum_i \lambda'_i O_i P_i)(\textstyle\sum_j O_j P_j)^{-1} \tag{47}$$

We check the basis transformation by Rodrigues formula [3]:
$$\mathbb{P}_n(x) = \frac{1}{\omega} D^n[\omega B^n(x)]$$

Where $\omega$ defined by the relation $\frac{\omega}{\omega'} = \frac{A-B'}{B}$ with $A$ as a polynomial of first degree.

If we choose monomial $B^n(x)$ as the original basis of vector space:
$$[1, B(x), B^2(x), \ldots, B^n(x)]$$
The Rodrigues formula could be chosen as the action of operator:
$$O_n = \frac{1}{\omega} D^n[\omega \cdot] \tag{48}$$
On these basis. Therefor it is a special case of operators of separated basis transformation.

The suitable operator with eigenfunctions $B^n(x)$ can be presented as
$$\frac{B(x)}{B'(x)} D\, B^n(x) = n B^n(x) \tag{49}$$
Where $B'(x)$ denoted as the derivative of $B(x)$. Thus the term $\frac{B(x)}{B'(x)} D$ should replace $\mathfrak{D}$ in (35):
$$P_k = \prod_{l=1}^{N} \frac{\frac{B(x)}{B'(x)} D - \lambda_l}{\lambda_k - \lambda_l} \tag{50}$$

Respect to (48) by replacing $O_n$ by Rodrigues formula and $P_n$ by equation (50) we get:
$$\mathfrak{D}' = (\sum_i \lambda'_i \frac{1}{\omega} D^n [\omega \prod_{l=1}^N \frac{\frac{B(x)}{B'(x)} D - \lambda_l}{\lambda_k - \lambda_l}])(\sum_j O_j P_j)^{-1}$$
Let $O^{-1} = (\sum_j O_j P_j)^{-1}$ then we obtain:
$$\mathfrak{D}' = (\sum_i \lambda'_i \frac{1}{\omega} D^n [\omega \prod_{l=1}^N \frac{\frac{B(x)}{B'(x)} D - \lambda_l}{\lambda_k - \lambda_l}]) O^{-1} \tag{51}$$
Taking into account the 2-dimensional space, and using (48) and (50) we have:
$$\lambda'_0 = 0, \quad O_1 = \frac{1}{\omega} D[\omega \cdot], \quad P_1 = \frac{\frac{B(x)}{B'(x)} D}{\lambda'_1} = \frac{1}{\lambda'_1} \frac{B(x)}{B'(x)} D$$
Thus the (51) reads as:
$$\mathfrak{D}' = \lambda'_1 \frac{1}{\omega} D \left[\omega \frac{1}{\lambda'_1} \frac{B(x)}{B'(x)} D\right] O^{-1} = \frac{1}{\omega} D \left[\omega \frac{B(x)}{B'(x)} D\right] O^{-1} \tag{52}$$
$$\mathfrak{D}' = (\frac{\omega}{\omega'} \frac{B}{B'} D + \frac{B'^2 - B''B}{B'^2} D + \frac{B}{B'} D^2) O^{-1}$$
If we assume $\frac{\omega}{\omega'} = \frac{A-B'}{B}$ ( as a crucial assumption in Rodrigues formula) this equation reduces to:
$$\mathfrak{D}' = \left\{\left(\frac{A}{B'} - 1\right) D + D - \frac{B''B}{B'^2} D + \frac{B}{B'} D^2\right\} O^{-1}$$
$$\mathfrak{D}' = (\frac{A}{B'} D - \frac{B''B}{B'^2} D + \frac{B}{B'} D^2) O^{-1} \tag{53}$$
Acting both side on $\mathbb{P}_1$ as the second eigenfunction of $\mathfrak{D}'$, we have:
$$\mathfrak{D}' \mathbb{P}_1 = (\frac{A}{B'} D - \frac{B''B}{B'^2} D + \frac{B}{B'} D^2) O^{-1} \mathbb{P}_1 \tag{54}$$
The term $O^{-1} \mathbb{P}_1$ equals $B(x)$, thus:
$$\mathfrak{D}' \mathbb{P}_1 = (\frac{A}{B'} D - \frac{B''B}{B'^2} D + \frac{B}{B'} D^2) B(x) \tag{55}$$
$$\mathfrak{D}' \mathbb{P}_1 = \left(\frac{A}{B'} - \frac{B''B}{B'^2} + \frac{B}{B'} D\right) B' = A - \frac{B''B}{B'} + \frac{B}{B'} (B'' + B'D)$$
$$\mathfrak{D}' \mathbb{P}_1 = BD + A$$
$$\mathfrak{D}' D^{-1} D \mathbb{P}_1 = BD + A \tag{56}$$
The term $D\mathbb{P}_1$ will be a constant $\alpha$, thus:
$$\alpha \mathfrak{D}' D^{-1} = BD + A$$
Or:
$$\alpha \mathfrak{D}' = BD^2 + AD \tag{57}$$

This implies that Rodrigues formula gives the solutions (or eigenfunctions) of the differential operator $BD^2 + AD$ and related differential equation up to a constant coefficient $\alpha$. i.e.,
$$\mathfrak{D}'y = \beta y \qquad (58)$$
Therefor any differential equation in the form of equation (58) can be solved by Rodrigues formula as a special case of separated basis transformation.
The following examples clarify this technique for some polynomials.

## Example 3.1  Laguerre differential equation

Let we intend to find the differential equation which corresponds to a set of linearly independent polynomials in variable $x$. For example, we are given a few first Laguerre polynomials i.e., $(1, 1 - x, ...)$ and we know the operator that maps the standard basis $(1, x, x^2, ...)$ to Laguerre basis i.e., operator presented in (41).
We can recover the corresponding Laguerre differential equation (operator) via the formula:
$$\mathfrak{D}' = (\sum_i \lambda'_i O_i P_i)(\sum_j O_j P_j)^{-1} \qquad (59)$$
*Proof in 2 dimension (first 2 polynomials)*

We restrict calculation in 2-dimensional polynomial space with basis $(1, x)$. These polynomials are transformed by $O_i$ to Laguerre polynomials in the same dimension i.e., $(1, -x + 1)$. Thus, the corresponding operator $xD$ will be transformed to Laguerre differential operator by the equation (41).
Substitution of $O_i$ by equation (41) and taking $\lambda'_i$ as the eigenvalues of Laguerre differential equation in 2-dimensional space of polynomials and replacing projection operators $P_i$ for basis $(1, x)$ by equation (35) into equation (37) results in:

$$\mathfrak{D}' = (\sum_{i=0}^{1} \lambda'_i O_i P_i)(\sum_{j=0}^{1} O_j P_j)^{-1} \qquad (60)$$

We have $\lambda'_i = \lambda_i$ and $\lambda_0 = 0$, $\lambda_1 = 1$. Then equation (60) reduces to:

$$\mathfrak{D}' = O_1 P_1 (O_0 P_0 + O_1 P_1)^{-1} \qquad (61)$$

By equations (35), (37) and (41) we obtain:
$$O_0 = 1 \ , \ O_1 = D - 1 \ , \ P_0 = \prod_{l=0}^{1} \frac{\mathfrak{D} - \lambda_1}{\lambda_0 - \lambda_1} = \frac{xD - 1}{-1} \ , \ P_1 = \prod_{l=0}^{1} \frac{\mathfrak{D} - \lambda_0}{\lambda_1 - \lambda_0} = \frac{xD}{1}$$

Therefor we have:
$$O_1 P_1 = (D - 1)xD = D + xD^2 - xD$$

And: $\quad \mathfrak{D}' = O_1 P_1 (O_0 P_0 + O_1 P_1)^{-1} = (D + xD^2 - xD)(D + xD^2 - 2xD + 1)^{-1} \qquad (62)$

If we denote the $(D + xD^2 - xD)$ as $\mathbb{D}$, we can reduce the equation (62) as follows:

$$\mathfrak{D}' = O_1 P_1 (O_0 P_0 + O_1 P_1)^{-1} = \mathbb{D}(\mathbb{D} - \mathfrak{D} + 1)^{-1} \qquad (63)$$

The term
$$O_0 P_0 + O_1 P_1 = \mathbb{D} - \mathfrak{D} + 1 = \hat{O} \qquad (64)$$
is the linear operator which transforms the basis $(1, x)$ to Laguerre basis $(1, -x + 1)$ and vice versa. If we restrict the action of operators to 2-dimensional polynomial space. Therefor we have:
$$\hat{O}^{-1}(-x + 1) = x \qquad (65)$$

If we act both sides of (47) on basis $(-x + 1)$ we get as well:
$$\mathfrak{D}'(-x + 1) = \mathbb{D}\hat{O}^{-1}(-x + 1) = \mathbb{D}x = -\mathbb{D}(-x + 1)$$
Or briefly:
$$\mathfrak{D}'(-x + 1) = -\mathbb{D}(-x + 1) \tag{66}$$
(Note that $-\mathbb{D}.1 = 0$)

The equation (66) implies that the action of both operators $\mathfrak{D}'$ and $-\mathbb{D}$ on basis $(1, -x + 1)$ are identical and therefor the simplest form of operator $\mathfrak{D}'$ which its eigenfunctions are Laguerre polynomials and its related transformation operators are $O_i$, reads as:
$$\mathfrak{D}' = -\mathbb{D} = -(xD^2 - xD + D) \tag{67}$$
This is the exactly the Laguerre differential equation with positive eigenvalues, i.e.:
$$-(xD^2 - xD + D)y = ny \tag{68}$$
Action of this operator on the first basis i.e., "1" gives 0 as the first eigenvalue and therefor the required conditions for validity of this differential operator are met.

*Proof in 3 dimension (first 3 polynomials)*

In 3-dimension with basis $(1, -x + 1, \frac{1}{2}(x^2 - 4x + 2))$ of Laguerre polynomial and $(1, x, x^2)$ of original basis, considering eigenvalues $\lambda_0 = 0$, $\lambda_1 = 1$, $\lambda_2 = 2$, the $\mathfrak{D}'$ reads as:
$$\mathfrak{D}' = (\sum_i \lambda'_i O_i P_i)(\sum_j O_j P_j)^{-1} = (O_1 P_1 + 2O_2 P_2)(O_0 P_0 + O_1 P_1 + O_2 P_2)^{-1} \tag{69}$$
$$\mathfrak{D}' = (O_1 P_1 + 2O_2 P_2)\hat{O}^{-1}$$
Here $\hat{O}^{-1}$ denotes the last term in (69). Acting both side on basis $(-x + 1)$ results in:
$$\mathfrak{D}'(-x + 1) = (O_1 P_1 + 2O_2 P_2)\hat{O}^{-1}(-x + 1) \tag{70}$$
Respect to (65) and the identity $P_2 x = 0$ we have:
$$\mathfrak{D}'(-x + 1) = O_1 P_1 x \tag{71}$$
In this dimension $P_1$ can be find as:
$$P_1 = \prod_{l=0}^{2} \frac{\mathfrak{D} - \lambda_l}{\lambda_1 - \lambda_l} = \left(\frac{xD - \lambda_0}{1 - \lambda_0}\right)\left(\frac{xD - \lambda_2}{1 - \lambda_2}\right) = xD\left(\frac{xD - 2}{1 - 2}\right) = -xD(xD - 2) \tag{72}$$
Then (71) reads as:
$$\mathfrak{D}'(-x + 1) = -(D - 1)xD(xD - 2)x$$
$$\mathfrak{D}'(-x + 1) = -(D - 1)xD(-x)$$
$$\mathfrak{D}'(-x + 1) = -[(D - 1)xD](-x + 1)$$
$$\mathfrak{D}'(-x + 1) = -(xD^2 - xD + D)(-x + 1)$$
This proves:
$$\mathfrak{D}' = -(xD^2 - xD + D) \tag{73}$$
As the Laguerre differential operator.

**Example 3.2   Hermite differential equation**

The same technique could be applied to derive Hermite differential equation by the formula (37). Because all $O_n$ that transforms basis $(1, x, x^2, x^3, \dots)$ to Hermite polynomials are equal to $O$ as is shown in (41), after getting $P_k$ by (35) and substitute them in (37) we have:

$$\mathfrak{D}' = (\sum_i \lambda'_i OP_i)(\sum_j OP_j)^{-1}$$
$$= O(\sum_i \lambda'_i P_i)(O \sum_j P_j)^{-1}$$
$$= O(\sum_i \lambda'_i P_i)(\sum_j P_j)^{-1} O^{-1}$$

Respect to $\sum_j P_j = 1$ we get:

$$\mathfrak{D}' = O(\sum_i \lambda'_i P_i)O^{-1} \tag{74}$$

Expanding the sum for eigenvalues $\lambda'_i = \lambda_i = 0,1$ and substitution of $O$ by equation (42) we have:

$$\mathfrak{D}' = e^{\frac{-D^2}{2}}(P_1)e^{\frac{D^2}{2}} \tag{75}$$

From (35) we calculate $P_1$ as:

$$P_1 = \prod_{l=0}^{1} \frac{\mathfrak{D} - \lambda_0}{\lambda_1 - \lambda_0} = \frac{\mathfrak{D} - 0}{1 - 0} = \mathfrak{D} \tag{76}$$

We know $(1, x, x^2, x^3, \ldots)$ are the eigenfunctions of $xD$, thus by $\mathfrak{D} = xD$ we have:

$$\mathfrak{D}' = e^{\frac{-D^2}{2}}(xD)e^{\frac{D^2}{2}} \tag{77}$$

This equation can be interpreted as a similarity transformation that maps $xD$ into $\mathfrak{D}'$ after basis changes. This will be hold just for the cases that eigenvalues are common between $xD$ and $\mathfrak{D}'$ as we see in Hermite and Laguerre differential equations.

Expansion of $e^{\frac{-D^2}{2}}$ and $e^{\frac{D^2}{2}}$ results in:

$$\mathfrak{D}' = (1 - \frac{D^2}{2} + \frac{D^4}{8} + \cdots)(xD)(1 + \frac{D^2}{2} + \frac{D^4}{8} + \cdots) \tag{78}$$

$$\mathfrak{D}' = (1 - \frac{D^2}{2} + \frac{D^4}{8} + \cdots)(xD + x\frac{D^3}{2} + x\frac{D^5}{8} + \cdots)$$

$$\mathfrak{D}' = \left(xD - \frac{D^2}{2}xD + \frac{D^4}{8}xD + \cdots\right) + (1 - \frac{D^2}{2} + \frac{D^4}{8} + \cdots)(x\frac{D^3}{2} + x\frac{D^5}{8} + \cdots)$$

In the 2-dimensional space of polynomials the orders higher than 2 for $D^n$ will be omitted, as it could be verified by action of both side on basis $x$. By omitting the higher orders, we obtain:

$$\mathfrak{D}' = (1 - \frac{D^2}{2})xD$$

$$\mathfrak{D}' = xD - \frac{D^2}{2}xD = xD - \frac{1}{2}D(D + xD^2)$$

$$\mathfrak{D}' = xD - \frac{1}{2}(D^2 + D^2 + xD^3)$$

Omitting $xD^3$ results in:

$$\mathfrak{D}' = xD - D^2 = -(D^2 - xD) \tag{79}$$

This is the well-known Hermit probabilist's Hermite differential operator with Hermite polynomial as its eigenfunctions and positive eigenvalues 0, 1, 2, …. as its eigenvalues.

**Example 3.3  Legendre differential equation**

For Legendre polynomials we have:

$$\mathbb{P}_n(x) = \frac{1}{2^n n!} D^n (x^2 - 1)^n \tag{80}$$

That transforms the basis set $S = \{1, (x^2 - 1), (x^2 - 1)^2, \ldots\}$ to Legendre polynomials. We can choose the appropriate operator $\mathfrak{D}$ whose eigenfunctions are these basis. Simply we write:

$$\mathfrak{D} = \frac{x^2 - 1}{2x} D \tag{81}$$

Eigenfunctions of this operator are members of the set $S$.
The transforming operator is

$$O_n = \frac{1}{2^n n!} D^n \tag{82}$$

In this case the eigenvalues of $\mathfrak{D}$ and $\mathfrak{D}'$ (Legendre differential operator) are not identical and therefor the similarity transformation is not valid. However, we can apply the equation (37) after determining the $P_i$ from (35).

For calculating $P_i$ in 2-dimension, we have:

$$P_0 = \prod_{l=0}^{1} \frac{\mathfrak{D}-\lambda_1}{\lambda_0-\lambda_1} = \frac{\frac{x^2-1}{2x}D-1}{-1} = 1 - \frac{x^2-1}{2x}D \tag{83}$$

$$P_1 = \prod_{l=0}^{1} \frac{\mathfrak{D}-\lambda_0}{\lambda_1-\lambda_0} = \frac{\frac{x^2-1}{2x}D-0}{1-0} = \frac{x^2-1}{2x}D \tag{84}$$

Noe we get:

$$\sum_i \lambda'_i O_i P_i = 2O_1 \frac{x^2-1}{2x}D = 2\left(\frac{1}{2}D\frac{x^2-1}{2x}D\right) = D\frac{x^2-1}{2x}D \tag{85}$$

From equation (37) and (85) we get:

$$\mathfrak{D}' = (\sum_i \lambda'_i O_i P_i)(\sum_j O_j P_j)^{-1}$$

$$\mathfrak{D}' = D\frac{x^2-1}{2x}D(\sum_j O_j P_j)^{-1} = D\frac{x^2-1}{2x}DO^{-1} \tag{86}$$

Where $O^{-1} = (\sum_j O_j P_j)^{-1}$.

By action of both sides of (86) on basis $x$ as the second basis of Legendre polynomials in two dimension, we have:

$$\mathfrak{D}'x = D\frac{x^2-1}{2x}DO^{-1}x \tag{87}$$

Respect to $O^{-1}x = x^2 - 1$, (87) reads as:

$$\mathfrak{D}'x = D\frac{x^2-1}{2x}D(x^2-1)$$

$$\mathfrak{D}'(D^{-1}D)x = D\frac{x^2-1}{2x}(2x)$$

$$\mathfrak{D}'D^{-1}(Dx) = D(x^2-1)$$

$$\mathfrak{D}'D^{-1} = D(x^2-1)$$

$$\mathfrak{D}' = D(x^2-1)D = -D(1-x^2)D \tag{88}$$

Expansion of (88) gives:

$$\mathfrak{D}' = -[(1-x^2)D^2 - 2xD]$$

Which is the Legendre differential operator with positive eigenvalues $n(n+1)$.

### 3.3 Solutions to Differential Equations by Raising Operator Method

In this section we apply the raising operator $A^+$ that is defined by:

$$\mathbb{P}_{n+1} = A^+ \mathbb{P}_n \tag{89}$$

and separated basis transformation method to solve some special differential equations. First, we find the trivial raising operator for the set of basis $(1, B(x), B^2(x), ... B^n(x))$. The raising operator for these basis is $B(x)$, because:

$$B(x)B^m(x) = B^{m+1}(x) \tag{90}$$

The corresponding raising operator for basis $\mathbb{P}_n$, can be derived by similarity transformation as follows:

$$A^+ = OB(x)O^{-1} \tag{91}$$

Where $O$ is the operator of basis change from $B^n(x)$ to $\mathbb{P}_n$ by $\mathbb{P}_n = OB^n(x)$.
defined . We start with a known differential equation and first two trivial solutions i.e., the first two eigenfunctions $1, \mathbb{P}_0$ with the lowest eigenvalues. By the definition of raising operator $A^+$ defined by (80) we derive this operator by restriction to 2 dimension of polynomial space and using the first two terms of $\sum_j O_j P_j$ and Forbenius covariant operator to retrieve entire eigenfunction

(solutions) of the considered differential equation. This method is introduced through the following examples.

**Example 3.4**

As an example, for Laguerre differential equation, if we know the first two monomial i.e., $\mathbb{L}_0 = 1$ and $\mathbb{L}_1 = -x + 1$ as the trivial eigenfunctions, respect to equation (91) the raising operator is
$$A^+ = OB(x)O^{-1}$$
Where operator $O$ transforms the basis $[x^n]$ to Laguerre polynomials $\mathbb{L}_n$.
For Laguerre differential equation by $B = x$, the raising operator appears as
$$A^+ = OxO^{-1} = (\textstyle\sum_j O_j P_j)xO^{-1} \tag{92}$$
$$A^+ = (O_1 P_1)xO^{-1}$$
And acting both side on $\mathbf{1}$ as the first monomial we get
$$A^+ . \mathbf{1} = (O_1 P_1)xO^{-1}. \mathbf{1} \tag{93}$$
By $O^{-1}. \mathbf{1} = \mathbf{1}$ and $P_1 x = x$ and by the $O_1 = D - 1$, this equation yields
$$A^+ . \mathbf{1} = (D-1)x$$
$$\mathbf{1} = (A^+)^{-1}(D-1)x$$
The action of operator $(A^+)^{-1}(D-1)$ on $x$ is the same as $D$, then we have the identity:

$$(A^+)^{-1}(D-1) = D \tag{94}$$
$$(A^+)^{-1}(D-1)D^{-1} = \mathbf{1}$$
$$(A^+)^{-1}(1 - D^{-1}) = \mathbf{1}$$

And this gives,
$$A^+ = 1 - D^{-1} \tag{95}$$
Applying this operator on the first two Laguerre polynomials gives the nth solution

$$\mathbb{L}_n = (A^+)^n . \mathbf{1} = (1 - D^{-1})^n . \mathbf{1} \tag{96}$$
This method can be applied for any differential operator to find its eigenfunctions or ordered solutions.

**Example 3.5**

For Hermite differential equation to derive $O_1$ due to equation (37) for 2 dimension we have

$$\mathfrak{D}'_H = (\lambda'_0 O_0 P_0 + \lambda'_1 O_1 P_1)(\textstyle\sum_j OP_j)^{-1}$$
We know $\lambda'_0 = 0$, $\lambda'_1 = 1$ and $O = (\textstyle\sum_j OP_j)^{-1}$, therefore:

$$\mathfrak{D}'_H = O_1 P_1 O^{-1} \tag{97}$$
Acting both side on second basis $x$ with the definition for projection operator $P_1$ from equation (76), gives

$$\mathfrak{D}'_H x = O_1 P_1 O^{-1} x$$
$$\mathfrak{D}'_H x = O_1 P_1 x$$
$$\mathfrak{D}'_H x = O_1 (xD)x$$
This equation shows both operators in the equation are equivalent:

$$\mathfrak{D}'_H = O_1(xD) \tag{98}$$

Substitution for $\mathfrak{D}'_H$ and action of $D^{-1}$ on both sides, yields

$$(xD - D^2)D^{-1} = O_1(xD)D^{-1}$$
$$x - D = O_1 x$$

Or
$$(x - D)x^{-1} = O_1 \tag{99}$$

Respect to $A^+ = O_1 B O^{-1}$ we get:

$$A^+ = (x - D)x^{-1} x O^{-1}$$

Acting both side on $\mathbf{1}$ as the first $\mathbf{1}$ basis

$$A^+ . \mathbf{1} = (x - D)x^{-1} x O^{-1} . \mathbf{1} \tag{100}$$
$$A^+ . \mathbf{1} = (x - D)x^{-1} x . \mathbf{1}$$

Thus, we have

$$A^+ = x - D$$

With raising operator, we derive all Hermits eigenfunctions as solutions to its differential equation:

$$\mathbb{H}_{e_n} = (A^+)^n . \mathbf{1} = (x - D)^n . \mathbf{1} \tag{101}$$

## Conclusion

We introduce a method of basis vector transformation of a vector space that applies separated operators to change the basis. Based on the method, a new equation for transforming the projection operators from a basis to another basis is presented. If the basis vectors considered as the eigenfunction $\varphi_n(x)$ of some differential operator $\mathfrak{D}$ in polynomial space that are separately transformed by operators $O_n$ to $\mathbb{P}_n(x)$, then by applying this method with the Forbenius covariants presented as projection operators, we derive the corresponding differential operator $\mathfrak{D}'$ for $\mathbb{P}_n$ as its eigenfunctions. This result shows that the Rodrigues formula is a special case of separated basis transformation method. The authenticity of the method has been verified by examples of classical orthogonal polynomials such as Hermite, Legendre and Laguerre polynomials. and substitution of them into the  that the  applicability of this method. A solution to some important differential equations by means of this method and raising operators in polynomial space has been introduced. It is expected that separated basis transformation could be applied in the theory of vector spaces and specifically in Hilbert space and quantum mechanics.

### Declarations

- The author confirms sole responsibility for the study conception and manuscript preparation.
- No funding was received for conducting this study.
- The author declares that he has no conflict of interest.

## References

1) Ismail, Mourad, Mourad EH Ismail, and Walter van Assche. *Classical and quantum orthogonal polynomials in one variable*. Vol. 13. Cambridge university press, 2005.
2) Arfken, George B., and Hans-Jurgen Weber. "Mathematical methods for physicists." (1972).
3) Rasala, Richard. "The Rodrigues formula and polynomial differential operators." *Journal of Mathematical Analysis and Applications* 84.2 (1981): 443-482.

4) Horn, Roger A., and Charles R. Johnson. "Topics in matrix analysis Cambridge university press." *Cambridge, UK* (1991). Pp(403)
5) Benzi, Michele, and Valeria Simoncini, eds. *Exploiting hidden structure in matrix computations: algorithms and applications*. Springer, 2016.
6) Doman, Brian George Spencer. *The classical orthogonal polynomials*. World Scientific, 2015.